\documentclass[twoside]{article}
\usepackage{amssymb, amsmath}
\begin{document}
\title{\vspace{-1in}\parbox{\linewidth}{\footnotesize\noindent
{\sc  Electronic Journal of Differential Equations},
Vol. {\bf 2002}(2002), No. ??, pp. 1--20. \newline
ISSN: 1072-6691. URL: http://ejde.math.swt.edu or http://ejde.math.unt.edu
\newline ftp  ejde.math.swt.edu  (login: ftp)}
 \vspace{\bigskipamount} \\
  The Vlasov equation with strong magnetic field and
  oscillating electric field as a model for isotop resonant separation
\thanks{ 
\indent
Pseudo-differential equations, Isotop separation.
\hfil\break\indent
\copyright 2002 Southwest Texas State University. \hfil\break\indent
Submitted November 24, 2001. Published ??.} }
\date{}
\author{Emmanuel Fr\'enod \& Fr\'ed\'erique Watbled}
\maketitle

\begin{abstract}
  We study the qualitative behavior of solutions to the Vlasov equation
  with strong external magnetic field and oscillating electric
  field. This model is relevant to the understanding of isotop
  resonant separation.  We show that the effective equation is
  a kinetic equation with a memory term.  This memory term involves
  a pseudo-differential operator whose kernel is characterized by an
  integral equation involving Bessel functions. The kernel is explicitly
  given in some particular cases.
\end{abstract}

\newcommand{\grad}{\mathop{\rm grad}}
\newtheorem{theorem}{Theorem}[section]
\newtheorem{remark}[theorem]{Remark}
\newtheorem{lemma}[theorem]{Lemma}

\renewcommand{\theequation}{\thesection.\arabic{equation}}
\catcode`@=11
\@addtoreset{equation}{section}
\catcode`@=12
\newtheorem{thm}{Theorem}[section]
\newtheorem{lem}[thm]{Lemma}
\newtheorem{rem}{Remark}[section]

\section{Introduction}

This paper gives a mathematical analysis of a model related to
 isotop resonant separation. We undertake this model using
homogenization methods applied to the Vlasov equation.

The mathematical framework of this paper is the investigation of the influence
of oscillations generated by strong fields in the Vlasov equation.
It completes the works led in
Fr\'enod \cite{frenod:1994},
Fr\'enod and Hamdache \cite{frenod/hamdache:1996},
Fr\'enod and Sonnendr\"ucker
\cite{frenod/sonnendrucker:PROC1998,frenod/sonnendrucker:1997,
frenod/sonnendrucker:CRAS99,frenod/sonnendrucker:1998,
frenod/sonnendrucker:1999},
Fr\'enod, Raviart and Sonnendr\"ucker \cite{FRS:1999},
Golse and Saint Raymond \cite{GSRcras,GSR1},
Saint Raymond \cite{saintraymond:2000},
Brenier \cite{brenier:2000},
Grenier \cite{grenier:1997,grenier2:1997},
Jabin \cite{jabin:2000},
Schochet \cite{schochet:1994},
Joly, M\'etivier and Rauch \cite{JMRDDD}.
We also refer to mathematical and physical works where similar
methods are used: \cite{gues93, JMR93, JMR96, serre91,
littlejohn:1981, lee:1983, dubin/etal:1983, cohen:1985,lochak/meunier:1988,grad:1971}.
The goal here is to exhibit the effect of the interaction of the
oscillations induced
by a strong magnetic field with the oscillations of the electric field.
In the case when both oscillation frequencies are the same, resonant phenomena
appear leading to memory effects in the effective equation. Concerning this
topic of memory effects, this work is an important step in the understanding
of non local homogenization previously analyzed in Sanchez-Palencia
\cite{sanchez-palencia:1978}, Tartar \cite{tartar:1989,tartar:1990},
Lions \cite{lions:HNL},
Amirat, Hamdache and Ziani \cite{AHZ:1989,AHZ:1992} and
Alexandre \cite{alexandre}, in the sense that we are able, here, to exhibit
explicit memory terms for a physically relevant problem.

 From the physical point of view, this study contributes to the understanding
of phenomena appearing during isotop resonant separation experiments.
Isotop resonant separation consists, in a plasma made of several kinds of ions,
in heating only the ions of one given species.
Those ions can be then easily extracted from the plasma by ad hoc devices.
In order to reach this goal, the plasma is,
inside a cylinder (with length about 1 m
and radius about 10 cm), submitted to a strong,
static and homogeneous magnetic field $B$ (about 0,4 T). Under its action,
each particle moves helicoidally around the magnetic lines with pulsation
$$ \omega_c=\frac {q|B|} {m}\,,
$$
called cyclotron pulsation, where $q$ and $m$ stand for electric charge
and mass of the considered particle.
Let $\tilde m$ and $\tilde q$ be the mass and the charge of the ions to be
heated. If an electric field $E$, oscillating with pulsation
$\frac{\tilde q |B|} {\tilde m}$ is moreover applied to the plasma,
the particles of the considered isotop are resonating with it and then
acquire energy.
For a detailed description of isotop resonant separation, we refer to
Omn\`es \cite{omnes:1999}, Louvet and Omn\`es \cite{louvet/omnes},
Schmitt \cite{schmitt:1973},
Dawson {\it et al.} \cite{dawson},
Compant La Fontaine and Louvet \cite{compant1,compant2}.

As a model enabling us to understand some aspects of isotop
resonant separation, we introduce a small parameter $\varepsilon$
and  consider the Vlasov equation
\begin{equation}\label{equ1}
\begin{gathered}
\partial_t f^\varepsilon + v\cdot\nabla_x f^\varepsilon
+ \big(E(t,\frac{t}{\varepsilon},x)+
v\times \big(B(t,\frac{t}{\varepsilon},x)+\frac{\mathcal{M}}{\varepsilon}\big)\big)
\cdot\nabla_v f^\varepsilon =0,\\
 f^\varepsilon(0,x,v)=f_0(x,v),\\
\end{gathered}
\end{equation}
where $\frac{\mathcal{M}}{\varepsilon}$ is a strong magnetic field,
$\mathcal{M}=e_1$ denoting the first
vector of the canonical basis $(e_1,e_2,e_3)$ of $\mathbb{R}^3$,
and where
$E(t,\frac{t}{\varepsilon},x)$ and $B(t,\frac{t}{\varepsilon},x)$
are fast oscillating electric and magnetic fields.
In this equation  $f^\varepsilon\equiv f^\varepsilon (t,x,v)$;
$t\in [0,T[$, where $T<+\infty$, is the time,
$x=(x_1,x_2,x_3)$ stands for the position and $v=(v_1,v_2,v_3)$ is the
velocity.
This equation models the evolution of the ions to be heated without
taking into account their interactions and the interactions
with the other ones. This is not completely unreasonable
since the ion density met in isotop separation experiment
is relatively low (about $10^{10}$ ions.cm$^{-1}$).
Yet, a model taking into account self induced forces
(like Vlasov-Poisson) will be considered in a forthcoming
paper.

We introduce the notations $\Omega=\mathbb{R}^3_x\times \mathbb{R}^3_v$
and ${\cal Q}=[0,T[\times\Omega$.
The initial data satisfies
\begin{equation}\label{equ2}
 f_0\geq 0,\quad 0<\int_\Omega f_0^2 \,dx \,dv <\infty.
\end{equation}
The electric and magnetic fields $E(t,\tau,x)$,
$B(t,\tau,x)$, are $C^\infty$ and
$2\pi$ periodic in $\tau$.
Under these assumptions we have the classical a priori estimate
that there exists a constant $C$ independent of $\varepsilon$ such that the
solution $f^\varepsilon$ of the Vlasov equation (\ref{equ1})
satisfies
$$\|f^\varepsilon\|_{L^\infty(0,T;L^2(\Omega))}\leq C.
$$
We deduce that up to a subsequence still denoted by $\varepsilon$,
$$f^\varepsilon \rightharpoonup  f \quad\mbox{in }
L^\infty(0,T;L^2(\Omega)) \quad\mbox{weak-}\star.
$$
The aim of this paper is to find out
an equation satisfied by the limit $f$.
Let us introduce some more notations:
for any vector $V$ in $\mathbb{R}^3$ we denote by $V_{\parallel}=V_1 e_1$,
$V_\perp=V_2 e_2+V_3 e_3$, respectively, the parallel and perpendicular
components of $V$ with respect to $e_1$.
If $\tau$ is in $[0,2\pi]$, we denote by $r(V,\tau)$ the image of $V$ under
the rotation of angle $\tau$ around the axis $\mathcal{M}=e_1$ and we define:
$$\tilde E (t,x)=\frac{1}{2\pi}\int_0^{2\pi} r\big(E(t,\tau,x),\tau\big)\,d\tau,$$
$$\tilde B (t,x)=\frac{1}{2\pi}\int_0^{2\pi} r\big(B(t,\tau,x),\tau\big)\,d\tau.$$
We shall denote by $J_0$ (respectively $J_1$) the Bessel function of order
zero (respectively of order one):
\begin{gather*}
J_0(z)=\frac{1}{\pi}\int_0^\pi \cos(z\sin\tau)\,d\tau
=\frac{1}{\pi}\int_0^\pi \cos(-z\cos\tau)\,d\tau,\\
J_1(z)=\frac{1}{\pi}\int_0^\pi \cos(z\sin\tau-\tau)\,d\tau
=-J_0'(z). \end{gather*}
 About definitions and properties of
Bessel functions, see for instance \cite{lebedev} or
\cite{watson}. \break To finish with the notations let us precise
that we shall use the notation of Taylor (\cite{taylor}) for
pseudo\-differential operators:
$$K(z;D)h(z)=\frac{1}{(2\pi)^{3/2}}
\int_{\mathbb{R}^3}K(z,ik)\mathcal{F} h(k)e^{iz\cdot k}\,dk.$$ Now
we are ready to state our  main result:
\begin{thm}[Main result] \label{th1}
We assume that $f_0(x,v)$ depends only on $x$,
$v_\parallel$, and the modulus of $v_\perp$, and satisfies
(\ref{equ2}). We assume also that $\tilde E(t,x)$ does not depend
on $x_\parallel$, that $\tilde E_\parallel=0$ and that $\tilde
B=0$. Then the sequence $(f^\varepsilon)$ of solutions of
(\ref{equ1}) satisfies, for any $T>0$,
$$f^\varepsilon \rightharpoonup  f \quad\mbox{in }
L^\infty(0,T;L^2(\Omega)) \quad\mbox{weak-}\star,
$$
and $f$ formally satisfies the following partial integro
differential equation
\begin{equation}\label{equ3}
\begin{gathered}
 \partial_t f(t,x,v) + v_1\cdot\partial_{x_1} f(t,x,v)
=\int_0^t K(s,t,y;D)f(s,x,v)~ds, \\
 f(0,x,v)=f_0(x,v),\\
\end{gathered}
\end{equation}
where $K(s,t,y;D)$ is a pseudo-differential operator parametrized
by $(s,t)$ in $[0,T[^2$ and $y=(v_1,x_2,x_3)$ in $\mathbb{R}^3$,
acting on functions of $z=(x_1,v_2,v_3)$ variables. The kernel
\begin{equation}\label{DefK}
K(s,t,y;ik) = \exp(-ik_1y_1(t-s)) \tilde K(s,t,y;ik)
\end{equation}
where $\tilde K(s,t,y;ik)$ is the unique solution of the Volterra
equation
\begin{multline}\label{EqKtilde}
 \int_s^t J_0\big(|k_\perp||L(s,\sigma,y)|\big)
\tilde K(\sigma,t,y,ik)\,d\sigma \\
 =-\frac{L(s,t,y)\cdot\tilde E(t,y)}{|L(s,t,y)|}|k_\perp|
J_1\big(|k_\perp||L(s,t,y)|\big)
\end{multline} %\end{equation}
with
\begin{equation}\label{DefL}
L(s,t,y)=\int_s^t \tilde{E}(\sigma,y)d\sigma
=\int_s^t \tilde{E}_\perp(\sigma,y)d\sigma.
\end{equation}
\end{thm}

\begin{rem} \rm
The assumption on variable-dependancy as the one on $\tilde E_\parallel$
we make in the heading of the Theorem is to ensure a transversality
property required in order to apply non local homogenization
methods.
Nevertheless, from physical point of view, as original repartition
functions are often maxwellian, considering that $f_0$ depends only
on $x$, $v_\parallel$ and the modulus of $v_\perp$ is relevant.
The two facts $\tilde E_\parallel=0$ and $\tilde B=0$ can be
easily realized experimentally, of course as soon as the self induced fields
are neglected.
\end{rem}

\begin{rem} \rm
Equation (\ref{equ3}) does not contain oscillations. It
only contains the mean effect, via the memory term,
of the oscillations contained in equation (\ref{equ1}).
According to our knowledge, it is the first time that such a mean
model is exhibited to describe the dynamic of ions during isotop resonant
separation.
\end{rem}

In some cases, we may obtain explicit expressions of the kernel $K$:

\begin{thm}\label{Thsuppl1}
In the particular case where
$$E_\perp(t,\tau,x_\perp)={\cal E}(t,x_\perp)g(\tau,x_\perp),$$
with ${\cal E}$ $\mathbb{R}$-valued and under the assumptions of Theorem \ref{th1},
the kernel is exactly equal to
\begin{equation}\label{NNNNN}
\begin{aligned}
 K&(s,t,y,ik)\\
 &=-\exp(-ik_1y_1(t-s))\,
{\cal E}(t,y){\cal E}(s,y)\; |k_\perp||\tilde g(y)|
\frac{J_1\big(|k_\perp||\tilde g(y)|\int_s^t {\cal E}
(\sigma,y)\,d\sigma\big)} {\int_s^t {\cal E}(\sigma,y)\,d\sigma},
\end{aligned}
\end{equation}
where
$$\tilde g(y)=\frac{1}{2\pi}\int_0^{2\pi} r\big(g(\tau,x),\tau\big)\,d\tau.$$
\end{thm}

\begin{rem} \rm
As it is stated in Alexandre \cite{alexandre},
it is immediate to check that the kernel $K$ appearing in (\ref{NNNNN})
satisfies:
$$|K(s,t,y,ik)|\leq \|{\cal E}\|_\infty |\tilde g(y)|^2(1+|k|^2).
$$
Moreover, we notice that for $s\neq t,$ the kernel $K(s,t,y,ik)$ decreases
like $\sqrt{|k|}$ when $k$ goes to infinity.
\end{rem}

The following Theorem is a simplified version of the previous one,
but as its proof is much simpler; we state it separately.

\begin{thm}\label{Thsuppl2}
In the very particular case of an electric field independent of time:
\begin{equation}
E_\perp = E_\perp(\tau,x_\perp)
\end{equation}
and under the assumptions of Theorem \ref{th1}, we have:
\begin{equation}\label{NNNNN2}
K(s,t,y,ik)=-\exp(-ik_1y_1(t-s))\, |k_\perp| |\tilde E_\perp(y)|
\frac{J_1\big(|k_\perp| |\tilde E_\perp(y)|(t-s) \big) }{t-s}.
\end{equation}
\end{thm}
The proof uses the notion of two scale convergence introduced
by N'Guetseng \cite{nguetseng:1989} and Allaire \cite{allaire:1992}.
Their result is the following:
\begin{thm}[N'Guetseng and Allaire] \label{th2}
If a sequence $(f^\varepsilon)$ is bounded in
\break
$L^\infty(0,T;L^2(\Omega))$,
then there exists a $2\pi$-periodic in $\tau$ profile
$F(t,\tau,x,v)$ in
\break
$L^\infty(0,T;L^\infty(\mathbb{R}_\tau;L^2(\Omega)))$
such that, for every $\psi(t,\tau,x,v)$ regular, compactly supported
with respect to $(t,x,v)$, and $2\pi$-periodic with respect to $\tau$,
we have, up to a subsequence,
$$\int_{\cal Q} f^\varepsilon \psi^\varepsilon ~dt ~dx ~dv
\rightarrow \int_{\cal Q} \int_0^{2\pi} F\psi \,d\tau ~dt ~dx
~dv,$$ where
$\psi^\varepsilon(t,x,v)=\psi(t,\frac{t}{\varepsilon},x,v)$. The
profile $F$ is called the $2\pi$-periodic two scale limit of
$f^\varepsilon$ and the link between $F$ and the weak-$\star$
limit $f$ is given by
$$\int_0^{2\pi} F(t,\tau,x,v)\,d\tau=f(t,x,v).$$

\end{thm}
It has been used by Fr\'enod and Sonnendr\"ucker
\cite{frenod/sonnendrucker:1997,frenod/sonnendrucker:1998,frenod/sonnendrucker:1999}
in the context of homogenization of the Vlasov equation but
under an assumption of strong convergence of the electric field.
Using the same ideas we obtain first the following result
concerning the two-scale limit of $f^\varepsilon$:
\begin{thm}\label{th3}
The sequence $(f^\varepsilon)$ of solutions of
(\ref{equ1}) two scale converges towards the $2\pi$-periodic in $\tau$
profile $F$ which is the unique solution
of
\begin{equation}\label{equ4}
\begin{gathered}
\partial_\tau F + (v\times \mathcal{M})\cdot\nabla_v F=0,\\
\partial_t F  + v_1\cdot\partial_{x_1}F
+ \big(r(\tilde{E},-\tau)+v\times r(\tilde{B},-\tau)\big)
\cdot\nabla_v F=0\\
 F(0,\tau,x,v)=\frac{1}{2\pi}f_0\big(x,r(v,\tau)\big).\\
\end{gathered}
\end{equation}
\end{thm}
Then following Alexandre
\cite{alexandre}
we apply a Fourier transform $\mathcal{F}$
to obtain an ordinary differential equation
satisfied by $\mathcal{F} F$.
Using the fact that
$$\mathcal{F} f(t,x,v)=\int_0^{2\pi} \mathcal{F} F(t,\tau,x,v)\,d\tau,$$
we obtain an ordinary differential equation
satisfied by $\mathcal{F} f$.
Eventually we introduce a Volterra equation as in Tartar
\cite{tartar:1989,tartar:1990}
and
Alexandre \cite{alexandre} to obtain our result.

\section{Scaling and qualitative study}

We exhibit the important parameters playing a role when charged
particles are submitted to a strong magnetic field and in view of which
we provide the scaling leading to equation (\ref{equ1}).

Before any scaling the evolution of the repartition function $f(t,x,v)$
representing at each time $t$ the particle density standing in $x$
and moving with velocity $v$, is given by the Vlasov equation
\begin{equation}
\partial_t f + v\cdot \nabla_x f +\frac{q}{m}(E+v\times B)\cdot\nabla_v f=0.
\end{equation}
We define now characteristic scales:
$\overline{t}$ is the characteristic time,
$\overline{L}$ the characteristic length and
$\overline{v}$ the characteristic velocity;
and rescaled variables: $t'$, $x'$ and $v'$  by
$t=\overline{t}t'$,
$x=\overline{L}x'$,
$v=\overline{v}v'$.
We also define scaling factors for the fields
$\overline{E}$ and $\overline{B}$ and rescaled fields
$E'(t',x')$ and $B'(t',x')$ by
$\overline{E}E'(t',x')=E(\overline{t}t',\overline{L}x')$
and
$\overline{B}B'(t',x')=B(\overline{t}t',\overline{L}x')$.
Lastly, defining a scaling factor $\overline{f}$ for the repartition function,
without forgetting that $f$ is a density on the phase space,
we define the rescaled repartition function $f'$ by
$$\overline{f}f'(t',x',v')=
\overline{L}^3\overline{v}^3 f(\overline{t}t',\overline{L}x',\overline{v}v').$$
The new repartition function is the solution to
\begin{equation}\label{2.6}
\partial_{t'} f' +
\frac{\overline{v}\overline{t}}{\overline{L}}v'\cdot \nabla_{x'} f'
+\big(\frac{q\overline{E}\overline{t}}{m\overline{v}}E'(t',x')
+\frac{q\overline{B}\overline{t}}{m}
v'\times B'(t',x')\big)\cdot\nabla_{v'} f=0.
\end{equation}
Let us introduce two parameters having an important physical signification:
$\overline{\omega_c}=\frac{q\overline{B}}{m}$ is the characteristic
cyclotron pulsation and
$\overline{a_L}=\frac{\overline{v}}{\overline{\omega_c}}$
the characteristic Larmor radius. Looking at (\ref{2.6}) in view of these parameters,
we get
\begin{equation}\label{etlabete}
\partial_{t'} f' +
\overline{t}\overline{\omega_c}\frac{\overline{a_L}}{\overline{L}}
v'\cdot \nabla_{x'} f'
+\big(\overline{t}\overline{\omega_c}
\frac{\overline{E}}{\overline{v}\overline{B}}E'(t',x')
+\overline{t}\overline{\omega_c}v'\times B'(t',x')\big)\cdot\nabla_{v'} f'=0.
\end{equation}
Now, introducing the small parameter $\varepsilon$, we set
\begin{equation}
\frac{\overline{a_L}}{\overline{L}}=\varepsilon {\rm\ and\ }
\overline{\omega_c}\overline{t}=\frac{1}{\varepsilon},
\end{equation}
which means that the observation length scale and the observation time scale
are respectively large in front of the Larmor radius and the cyclotronic
period.
This regime is relevant to describe the global behaviour
of the considered particles.
We also assume that the electric force is much smaller than the magnetic one.
It reads:
\begin{equation}
\frac{\overline{E}}{\overline{v}\overline{B}}=\varepsilon.
\end{equation}
Now, in order to model
the fact that the particles are submitted to an oscillating electric field,
we set that $E'$ writes $E'\big(t',\frac{t'}{\varepsilon},x'\big)$.
Concerning the magnetic field we assume that it is made of a constant field
perturbed by an oscillating one, it gives
$B'=\mathcal{M}+\varepsilon B''\big(t',\frac{t'}{\varepsilon},x'\big)$.
Using those assumptions and removing subscripts $''$ and $'$, equation
(\ref{etlabete}) leads to (\ref{equ1}).

\begin{rem} \rm
The goal of this remark is to show that the considered scaling
corresponds to physical situations.
If the ions are potassium and the magnetic field magnitude is $1T$,
the characteristic cyclotron pulsation is about $10^5$
($=\overline{\omega_c}$).
In experimental situations it is realistic to consider that
the ions stay $10^{-2}$ s in the device (we make this value
as reference time $\overline{t}$).
Beyond this, considering that the thermic velocity (we choose
as characteristic velocity $\overline{v}$)
of the ions is something like $10^3 m.s^{-1}$, the size order
of their Larmor radius is the mm ($\overline{a_L}\sim 10^{-3}$m).
As the device is about one meter long, we have
\begin{equation}
\frac{\overline{a_L}}{\overline{L}}\sim 10^{-3},~~
\overline{\omega_c}\overline{t}=10^3.
\end{equation}
Yet if the magnitude of the electric field is about
$1 V.m^{-1}$ we also have
\begin{equation}
\frac{\overline{E}}{\overline{v}\overline{B}}=10^{-3}.
\end{equation}
\end{rem}

As stated in Theorem \ref{th1}, the effective behaviour
of equation (\ref{equ3}) involves memory effect.
In order to give a way to understand why, we propose
to begin by studying the following simple problem:
analysing the movement of a particle (with mass and electric charge equal
to $1$)
with no velocity component parallel to $\mathcal{M}=e_1$,
and submitted to the magnetic field $\mathcal{M}/\varepsilon$.
Under the mere action of this magnetic field, the particle rotates around
$\mathcal{M}$ with pulsation $1/\varepsilon$.
In other words, its velocity and position write
$$\begin{array}{lr}
\displaystyle
 V(t)=r(v_0,-\frac{t}{\varepsilon}), &
 X(t)=x_0+\varepsilon\big(
r(v_0,\frac{\pi}{2}-\frac{t}{\varepsilon})+v_0\times\mathcal{M}\big),\\
\end{array}$$
where $v_0$ and $x_0$ are velocity and position of the particle at $t=0$
(we assume $v_0\cdot \mathcal{M}=0$).
As $\varepsilon$ goes to $0$, $X(t)$ tends to $x_0$.
Yet, $V(t)$ drives the circle $|v|=|v_0|$ faster and faster.
Consequently we could say that {\it the particle
occupies the whole circle $|v|=|v_0|$
and has forgotten its initial direction $\frac{v_0}{|v_0|}$}.
If now an oscillating electric field, writing for instance
$$E^\varepsilon(t)=r(e_2,-\frac{t}{\varepsilon})
=\begin{pmatrix}
0 \\ \cos(t/\varepsilon)\\ -\sin(t/\varepsilon) \end{pmatrix}
$$
is applied to the particle in addition to the magnetic field, the movement
of the particle is no more a rotation, but a spiral around the magnetic field
$\mathcal{M}$.
Indeed, the considered electric field does not modify the angular velocity
of the particle but only the modulus
$|v|=\sqrt{v_2^2+v_3^2}$, since we have
$$V(t)=r(v_0,-\frac{t}{\varepsilon})+tE^\varepsilon(t).$$
Observe that
$$\frac{d}{dt}|V(t)|^2=2(v_0\cdot e_2 +t).$$
If $v_0\cdot e_2 \geq 0$ the modulus $|v|$ of the velocity
increases during each rotation; if $v_0\cdot e_2 <0$
the modulus decreases during the first rotations
(until the time $t=-v_0\cdot e_2$)
and then it increases.
Consequently, we could say here that
{\it the dynamic of the particle strongly depends on the initial
value $v_0$}.
Of course this dependance is kept when
$\varepsilon$ goes to $0$.
But as we just saw the particle {\it forgets its initial
value} as $\varepsilon$ goes to $0$.
Because of this contradiction
we need to keep additional information, which is contained
in the memory term taking place in equation (\ref{equ3}). \smallskip

To finish this qualitative study, we shall explain the result
contained in Theorem \ref{th3} using formal asymptotic
expansion. If we assume the following ansatz
of $f^\varepsilon$
\begin{equation}
f^\varepsilon(t,x,v)=F_0(t,\frac{t}{\varepsilon},x,v)
+\varepsilon F_1(t,\frac{t}{\varepsilon},x,v)+\ldots,
\end{equation}
if we insert this in (\ref{equ1}) and identify the terms
at each order we get, at order $0$:
\begin{equation}\label{I}
\partial_\tau F_0 + (v\times \mathcal{M})\cdot\nabla_v F_0=0,
\end{equation}
and
\begin{multline}\label{II}
 \partial_\tau F_1 + (v\times \mathcal{M})\cdot\nabla_v F_1 \\
 =  -\big(\partial_t F_0  + v\cdot\nabla_{x}F_0
+ \big(E(t,\tau,x)+v\times B(t,\tau,x)\big)
\cdot\nabla_v F_0\big).
\end{multline}
We then see that the first term $F_0$ is nothing but the two scale limit
$F$ and that equation (\ref{I})
is nothing but (\ref{equ4}a).
The second equation (\ref {equ4}b) is given as a compatibility condition
on $F_0=F$ in order that (\ref{II}) has solutions.

%%%%%%%%%%%%%%%%%%%%%%%%%%%%%%%%%%%%%%%%%%%%%%%%
%SECTION 3

\section{Equation satisfied by the 2-scale limit.}

In this section we follow the procedure of Fr\'enod and Sonnendr\"ucker
\cite{frenod/sonnendrucker:1997}
to prove Theorem \ref{th3}.
We only sketch the proof and refer the reader to
\cite{
frenod/sonnendrucker:1997} for details.

Let $(f^\varepsilon)$ be a sequence of solutions of (\ref{equ1}).
As the sequence is bounded in
$L^\infty(0,T;L^2(\Omega))$,
by Theorem \ref{th2} there exists a $2\pi$-periodic in $\tau$
profile
$F(t,\tau,x,v)$ in $L^\infty(0,T;L^\infty(\mathbb{R}_\tau;L^2(\Omega)))$
such that for every $\psi(t,\tau,x,v)$ regular, compactly supported
with respect to $(t,x,v)$ and $2\pi$-periodic with respect to $\tau$,
we have, up to a subsequence,
\begin{equation}\label{equ5}
\int_{\cal{Q}} f^\varepsilon \psi^\varepsilon ~dt ~dx ~dv
\rightarrow \int_{\cal{Q}}\int_0^{2\pi} F\psi \,d\tau ~dt ~dx ~dv
\end{equation}
The proof of Theorem \ref{th3} is led in three steps:

\noindent{\it Step 1:} First we use a weak formulation of (\ref{equ1}) in
$\mathcal{D} '({\cal Q})$ with functions
$\psi^\varepsilon(t,x,v)=\psi(t,\frac{t}{\varepsilon},x,v)$ where
$\psi$ is regular with compact support in $(t,x,v)$ and
$2\pi$-periodic in $\tau$, which writes:
\begin{multline}\label{equ6}
 \int_{\cal{Q}} f^\varepsilon
\big(\partial_t\psi^\varepsilon + v\cdot\nabla_x\psi^\varepsilon
+\big(E^\varepsilon+v\times(B^\varepsilon +\frac{\mathcal{M}}{\varepsilon})\big)
\cdot\nabla_v\psi^\varepsilon\big)\,dt\,dx\,dv  \\
  =-\int_\Omega f_0(x,v)\psi^\varepsilon(0,x,v)\,dx\,dv.
\end{multline}
Notice that
$\displaystyle
\partial_t(\psi^\varepsilon)=\big(\partial_t\psi+\frac{1}{\varepsilon}\partial_\tau\psi\big)
^\varepsilon$.
Multiply (\ref{equ6}) by $\varepsilon$, then let $\varepsilon$ tend to $0$ and
apply the two scale
convergence
(\ref{equ5}) to deduce that $F$ belongs to the kernel of the singular
perturbation appearing in (\ref{equ1}), in other words:
\begin{equation}\label{equ7}
\partial_\tau F + (v\times \mathcal{M})\cdot\nabla_v F=0
\quad\mbox{in } \mathcal{D} '(\mathbb{R}_\tau\times\mathbb{R}_v^3)
\end{equation}
for almost every $(x,t)$ in $[0,T[\times\mathbb{R}_x^3$.

\noindent{\it Step 2:}
Next we use the following lemma
(Lemma 2.3 of \cite{frenod/sonnendrucker:1997}), which characterizes
the kernel of the singular perturbation.
\begin{lem}\label{lem1}
A function $F(\tau,v)\in L^\infty(\mathbb{R}_\tau,L^2(\mathbb{R}_v^3))$
$2\pi$-periodic in $\tau$
satisfies
$$\partial_\tau F + (v\times \mathcal{M})\cdot\nabla_v F=0
\quad\mbox{in } \mathcal{D}
'(\mathbb{R}_\tau\times\mathbb{R}_v^3)$$ if and only if there
exists a function $G\in L^2(\mathbb{R}_u^3)$ such that
$F(\tau,v)=G(r(v,\tau))$.
\end{lem}

According to this lemma there exists
a function $G$ in
$L^\infty(0,T;L^2(\mathbb{R}_x^3\times\mathbb{R}_u^3))$
such that
\begin{equation}\label{equ8}
F(t,\tau,x,v)=G(t,x,r(v,\tau)).
\end{equation}

\noindent{\it Step 3:}
We denote here $\Omega'=\mathbb{R}_x^3\times\mathbb{R}_u^3$, ${\cal Q}'=[0,T[\times \Omega'$.
The goal of this step is to project equation (\ref{equ1})
on the orthogonal of the kernel we identified above
and to pass again to the limit.
In order to achieve this we build test functions belonging to the
kernel in the following way.
For every regular compactly supported $\varphi(t,x,u)$
we consider the $2\pi$-periodic in $\tau$ function
$\psi(t,\tau,x,v)=\varphi(t,x,r(v,\tau))$, which in view of
Lemma \ref{lem1} satisfies
$$\partial_\tau \psi  + (v\times \mathcal{M})\cdot\nabla_v \psi=0. $$
We take $\psi^\varepsilon(t,x,v)=\psi(t,\frac{t}{\varepsilon},x,v)$ in the weak formulation
of the Vlasov equation (\ref{equ6}) obtained in the first step
so that
\begin{multline*}
 \int_{\cal{Q}} f^\varepsilon
\big((\partial_t\psi)^\varepsilon + (v\cdot\nabla_x\psi)^\varepsilon
+\big((E+v\times B)
\cdot\nabla_v\psi\big)^\varepsilon\big)\,dt\,dx\,dv \\
 =-\int_\Omega f_0(x,v)\psi (0,0,x,v)\,dx\,dv.
\end{multline*}
We let $\varepsilon$ tend to $0$, use the two scale convergence (\ref{equ5})
and the equality (\ref{equ8}) to deduce that
\begin{align*}
\int_{\cal{Q}}\int_0^{2\pi} G\big(t,x,r(v,\tau)\big)
\big(\partial_t\varphi + v\cdot\nabla_x\varphi
+r\big(E+v\times B,\tau\big)
\cdot \nabla_u\varphi \big)\,d\tau~dt~dx~dv\\
=-\int_\Omega f_0(x,v)\varphi (0,x,v)~dx~dv.\\
\end{align*}
Now we make the change of variables $u=r(v,\tau)$ and perform the integration
with respect to $\tau$
over $[0,2\pi]$.
We obtain in this way the equation satisfied by $G$:
\begin{lem}\label{lem2}
The function $G(t,x,u)$ linked to the $2\pi$-periodic profile
$F$ by (\ref{equ8}) is the unique solution of
\begin{equation}\label{equ9}
\begin{gathered}
 \partial_t G + u_1\cdot\partial_{x_1}G
+(\tilde{E}+u\times \tilde B)\cdot\nabla_u G=0\\
 G(0,x,u)=\frac{1}{2\pi}f_0(x,u).
\end{gathered}
\end{equation}
\end{lem}

The uniqueness of the solution of (\ref{equ9}) enables us to
deduce that the whole sequence $f^\varepsilon$ two-scale converges
to $F$ and, because of the link between $F$ and $f$, weak-$\star$
converges to $f$ (recall that $\int_0^{2\pi}
F(t,\tau,x,v)\,d\tau=f(t,x,v)$).

To  prove Theorem \ref{th3} we rewrite the equation
(\ref{equ9}) in terms of $F$ using the equality
$F(t,\tau,x,v)=G(t,x,r(v,\tau))$
and thus obtain the equation satisfied by $F$
since
$$r\big(\nabla_v F(t,\tau,x,v),\tau\big)=\nabla_uG\big(t,x,r(v,\tau)\big)$$
and
$$\big(r(v,\tau)\times \tilde B\big)\cdot r(\nabla_vF,\tau)
=\big(v\times r(\tilde B,-\tau)\big)\cdot \nabla_v F.$$
The Theorem is then proved. \hfill$\Box$

\begin{rem} \rm
At this stage we can integrate in $\tau$ over $[0,2\pi]$ the equation
\begin{gather*}
\partial_t F  + v_1\cdot\partial_{x_1}F
+ \big(r(\tilde{E},-\tau)+v\times r(\tilde B,-\tau)\big)\cdot\nabla_v F=0\\
 F(0,\tau,x,v)=\frac{1}{2\pi}f_0\big(x,r(v,\tau)\big)\\
\end{gather*}
and use the equality
$f(t,x,v)=\int_0^{2\pi} F(t,\tau,x,v)\,d\tau$
to get the following equation which is satisfied by $f$:
\begin{multline*}
 \partial_t f  + v_1\cdot\partial_{x_1}f
+ \big(\tilde{E}_\parallel+v\times \tilde B_\parallel\big)\nabla_v f\\
+  \big(\tilde E_\perp + v\times \tilde B_\perp\big)
\cdot \nabla_v a
+ \big(\tilde E\times\mathcal{M} + v\times (\tilde B\times\mathcal{M})\big)
\cdot \nabla_v b
=0,
\end{multline*}
where
$a(t,x,v)=\int_0^{2\pi} F(t,\tau,x,v)\cos \tau \,d\tau$, \quad
$b(t,x,v)=\int_0^{2\pi} F(t,\tau,x,v)\sin \tau \,d\tau$.
\end{rem}

%%%%%%%%%%%%%%%%%%%%%%%%%%%%%%%%%%%%%%%%%%%%%
%SECTION 4

\section{Equation satisfied by the weak-* limit.}

In this section we prove our main results Theorems
\ref{th1}, \ref{Thsuppl1} and \ref{Thsuppl2}
by using Fourier transform.
Recall that under the assumptions of Theorem \ref{th1}
the $2$-scale limit
$F$  satisfies the equation

\begin{equation}\label{equ10}
\begin{gathered}
\partial_t F  + v_1\cdot\partial_{x_1}F
+ r(\tilde{E},-\tau) \cdot \nabla_v F=0\\
 F(0,\tau,x,v)=\frac{1}{2\pi}f_0(x;v_\parallel,|v_\perp|),\\
\end{gathered}
\end{equation}
with
$$r(\tilde{E}(t,x_\perp),-\tau)=\left|\begin{array}{l}
0\\
\tilde{E}_2(t,x_\perp)\cos\tau +\tilde{E}_3(t,x_\perp)\sin\tau\\
-\tilde{E}_2(t,x_\perp)\sin\tau+\tilde{E}_3(t,x_\perp)\cos\tau\\
\end{array}\right.$$
Thanks to the hypothesis the only derivatives of $F$
involved in the equation are with respect to
$x_1$, $v_2$, $v_3$ and the coefficients only depend on
$v_1$, $x_2$, $x_3$. The transversality assumption required for non
local homogenization methods is realized.
For convenience we rename the variables by
$$\begin{array}{lr}
z=\left|\begin{array}{l}x_1\\v_2\\v_3\\\end{array}\right.,&
y=\left|\begin{array}{l}v_1\\x_2\\x_3\\\end{array}\right..
\end{array}$$
and define
\begin{gather*}
H(t,\tau,z,y)=F(t,\tau,x,v),\\
h(t,z,y)=f(t,x,v),\\
h_0(z,y)=h(0,z,y)=f(0,x,v)=f_0(x,v),\\
a(t,\tau,y)=\left|\begin{array}{l}
y_1\\
\tilde{E}_2 \cos\tau+\tilde{E}_3\sin\tau\\
-\tilde{E}_2\sin\tau+\tilde{E}_3\cos\tau\\
\end{array}\right.,\\
\end{gather*}
so that (\ref{equ10}) becomes
\begin{equation}\label{equ11}
\begin{gathered}
 \partial_t H(t,\tau,z,y)+a(t,\tau,y)
\cdot\nabla_z H(t,\tau,z,y)=0,\\
 H(0,\tau,z,y)=\frac{1}{2\pi}h_0(z,y).
\end{gathered}
\end{equation}
Applying a Fourier transform in the $z$ variable we get
the ordinary differential equation
\begin{equation}\label{equ12}
\begin{gathered}
 \partial_t \mathcal{F} H(t,\tau,k,y)
+ik\cdot a(t,\tau,y)\mathcal{F} H(t,\tau,k,y)=0,\\
 \mathcal{F} H(0,\tau,k,y)=\frac{1}{2\pi}\mathcal{F} h_0(k,y),
\end{gathered}
\end{equation}
which has the explicit solution
\begin{align*}
\mathcal{F} H(t,\tau,k,y)&=
\frac{1}{2\pi}\mathcal{F} h_0(k,y)
\exp\big(-ik\cdot\int_0^t a(\sigma,\tau,y)\,d\sigma\big) \\
&=\frac{1}{2\pi}\mathcal{F} h_0(k,y)\exp(-ik_1y_1t)
\exp\big(-ik\cdot\int_0^t r\big(\tilde
E(\sigma,y),-\tau\big)\,d\sigma\big).
\end{align*}
As
$$f(t,x,v)=\int_0^{2\pi} F(t,\tau,x,v)\,d\tau$$
we know that
$$\mathcal{F} h (t,k,y)=\int_0^{2\pi} \mathcal{F} H(t,\tau,k,y)\,d\tau,$$
hence
$$\mathcal{F} h(t,k,y)=
\mathcal{F} h_0(k,y)\exp(-ik_1y_1t) \frac{1}{2\pi}\int_0^{2\pi}
\exp\big(-ik\cdot \int_0^t r\big(\tilde
E(\sigma,y),-\tau\big)d\sigma\big)\,d\tau.$$
We introduce the quantities
\begin{gather}
L(s,t,y)=\int_s^t \tilde{E}(\sigma,y)d\sigma, \\
A(s,t,y,Z)=\frac{1}{2\pi}\int_0^{2\pi} \exp\big(-Z\cdot \int_s^t
r\big(\tilde E(\sigma,y),-\tau\big)\,d\sigma \big)\,d\tau,
\label{EqDefA}
\end{gather}
so that
$$\mathcal{F} h(t,k,y)=\mathcal{F} h_0(k,y)\exp(-ik_1y_1t)A(0,t,y,ik),$$
and differentiating with respect to $t$, the following equation
satisfied by $\mathcal{F} h$:
\begin{equation}\label{equ13}
\begin{gathered}
\partial_t\mathcal{F} h(t,k,y)
+ik_1y_1\mathcal{F} h(t,k,y)=\mathcal{F} h_0(k,y)\exp(-ik_1y_1t)\partial_t A(0,t,y,ik),\\
\mathcal{F} h(0,k,y)=\mathcal{F} h_0(k,y).
\end{gathered}
\end{equation}
Now we use the method of Tartar
\cite{tartar:1989,tartar:1990} (see also Alexandre \cite{alexandre})
to form an integro-differential operator.
We define
$$C(s,t,y,Z)=-\partial_t A(s,t,y;Z_2,Z_3)\exp(-Z_1y_1(t-s))
\quad\mbox{for }Z\in\mathbb{C}^3,$$
so that
\begin{equation}\label{equ14}
\partial_t\mathcal{F} h(t,k,y)
+ik_1y_1\mathcal{F} h(t,k,y)=-\mathcal{F} h_0(k,y)C(0,t,y,ik).
\end{equation}
We denote by
$D(s,t,y,Z)$ the solution of the Volterra-Green equation
\begin{equation}\label{equ15}
D(s,t,y,Z)-\int_s^t
C(s,\sigma,y,Z)D(\sigma,t,y,Z)\,d\sigma=C(s,t,y,Z).
\end{equation}
We replace $C(0,t,y,ik)$ in equation (\ref{equ14})
by its integral form, which gives
\begin{align*}
\partial_t\mathcal{F}& h(t,k,y)
+ik_1y_1\mathcal{F} h(t,k,y)\\
=&-\mathcal{F} h_0(k,y)D(0,t,y,ik)
+\int_0^t \mathcal{F} h_0(k,y)C(0,\sigma,y,ik)D(\sigma,t,y,ik)\,d\sigma\\
=&-\mathcal{F} h_0(k,y)D(0,t,y,ik)
-\int_0^t \partial_t \mathcal{F} h(\sigma,k,y)D(\sigma,t,y,ik)\,d\sigma\\
&-ik_1y_1 \int_0^t \mathcal{F} h(\sigma,k,y)D(\sigma,t,y,ik)\,d\sigma.
\end{align*}
Integrating by parts in the second term, observe
that $$A(t,t,y,Z)=1, D(t,t,y,Z)=-\partial_tA(t,t,y,Z)=0,$$
and we obtain eventually that
\begin{equation}\label{equ16}
\begin{gathered}
 \partial_t\mathcal{F} h(t,k,y)
+ik_1y_1\mathcal{F} h(t,k,y)=\int_0^t K(\sigma,t,y,ik)\mathcal{F} h(\sigma,k,y)\,d\sigma,\\
 \mathcal{F} h(0,k,y)=\mathcal{F} h_0(k,y),\\
\end{gathered}
\end{equation}
 where
$K(s,t,y,Z)= \partial_s D(s,t,y,Z)-ik_1y_1D(s,t,y,Z)$.

To simplify the expression of $K$ we define $\tilde D$ by
$$\tilde D(s,t,y,Z)=\exp\big(Z_1y_1(t-s)\big)D(s,t,y,Z).
$$
Then replacing $D$ in the expression of $K$
we obtain that
\begin{equation}\label{equ17}
K(s,t,y,Z)= \exp\big(-Z_1y_1(t-s)\big) \tilde K(s,t,y,Z),
\end{equation}
which is equation (\ref{DefK}), with
\begin{equation}\label{equ17bis}
\tilde K(s,t,y,Z) = \partial_s\tilde D(s,t,y,Z),
\end{equation}
and replacing $D$ in the Volterra equation (\ref{equ15})
we obtain that $\tilde D$ is the solution of the Volterra
equation
\begin{equation}\label{equ18}
\tilde D(s,t,y,Z)+\int_s^t \partial_t A(s,\sigma,y,Z) \tilde
D(\sigma,t,y,Z)\,d\sigma=-\partial_t A(s,t,y,Z).
\end{equation}
We observe that $\tilde D$ does not depend on $Z_1$.
Performing an integration by parts we get the following equation
satisfied by $\tilde K = \partial_s\tilde D$:
\begin{equation}\label{equ19}
\int_s^t A(s,\sigma,y,Z) \tilde K(\sigma,t,y,Z)\,d\sigma
=\partial_t A(s,t,y,Z).
\end{equation}

Let us now simplify the (\ref{EqDefA}) expression of $A$:
for fixed $s,t,y,k$, we choose $\alpha$ and $\beta$ in
$[0,2\pi]$ such that
\begin{gather*}
r(L(s,t,y),\alpha)=|L(s,t,y)|e_2,\\
k_\perp=|k_\perp|r(e_2,\beta).
\end{gather*}
Then
\begin{align*}
A(s,t,y,ik)&=
\frac{1}{2\pi}\int_0^{2\pi}
\exp\big(-ik\cdot
\int_s^t r\big(\tilde E(\sigma,y),-\tau\big)\,d\sigma \big)\,d\tau \\
&=\frac{1}{2\pi}\int_0^{2\pi}
\exp\big(-ik\cdot
 r\big(L(s,t,y),-\tau\big) \big)\,d\tau \\
&=\frac{1}{2\pi}\int_0^{2\pi}
\exp\big(-i|k_\perp||L(s,t,y)|
 r(e_2,\beta+\alpha+\tau)\cdot e_2\big)\,d\tau \\
&=\frac{1}{2\pi}\int_0^{2\pi}
\exp\big(-i|k_\perp||L(s,t,y)|
\cos\tau\big)\,d\tau \\
&=\frac{1}{2\pi}\int_0^{2\pi}
\cos\big(|k_\perp||L(s,t,y)|
\cos\tau\big)\,d\tau,
\end{align*}
so that
\begin{equation}\label{equ 25}
A(s,t,y,ik)=J_0\big(|k_\perp||L(s,t,y)|\big).
\end{equation}
Differentiating with respect to $t$, we obtain
\begin{equation}\label{equ 26}
\partial_tA(s,t,y,ik)=
-\frac{L(s,t,y)\cdot\tilde E(t,y)}{|L(s,t,y)|}|k_\perp|
J_1\big(|k_\perp||L(s,t,y)|\big).
\end{equation}
Hence $\tilde K(s,t,y,ik)$
is the unique solution of the equation
\begin{multline}\label{equ 27}
 \int_s^t J_0\big(|k_\perp||L(s,\sigma,y)|\big)
\tilde K(\sigma,t,y,ik)\,d\sigma \\
= -\frac{L(s,t,y)\cdot\tilde E(t,y)}{|L(s,t,y)|}|k_\perp|
J_1\big(|k_\perp||L(s,t,y)|\big)
\end{multline}
for every $s,t\in [0,T]$, that is (\ref{EqKtilde}).

Applying formally the inverse Fourier transform to (\ref{equ16}) we find that
$$\partial_t h(t,z,y)+y_1\cdot\partial_{z_1}h(t,z,y)
=\int_0^t K(\sigma,t,y;D)h(\sigma,z,y)\,d\sigma.$$
Now
$$h(t,z,y)=\int_0^{2\pi} H(t,\tau,z,y)\,d\tau
=\int_0^{2\pi} F(t,\tau,x,v)\,d\tau= f(t,x,v),$$ from which we
deduce (\ref{equ3}) satisfied by $f$, and completes the
proof of Theorem~\ref{th1}. \hfill$\Box$ \smallskip

We shall now treat the particular case of Theorem \ref{Thsuppl2}
where $E$ is independent of time:
in this case $\tilde E=\tilde E_\perp$ is independant of time too
and we have
$$L(s,t,y)=(t-s)\tilde E_\perp(y),$$
so that using the parity of $J_0$ and the imparity of
$J_1$ the last equation becomes
\begin{multline*}
 \int_s^t J_0\big(|k_\perp||\tilde E_\perp(y)|(\sigma-s)\big)
\tilde K(\sigma,t,y,ik)\,d\sigma \\
 =-|\tilde E_\perp(y)||k_\perp|
J_1\big(|k_\perp||\tilde E(y)|(t-s)\big)\quad\mbox{for every }s,t\in [0,T].
\end{multline*}
For sake of clarity we fix for the moment $y$ and $k$ and we set
$$C=|k_\perp||\tilde E_\perp(y)|.
$$
It is easy to see that for every $h\in]-T,+T[$, the function
$(s,t)\mapsto \tilde K(s+h,t+h)$ is solution of the
equation on $[-h,T-h]\times [-h,T-h]$, so that
$\tilde K(s+h,t+h)=\tilde K(s,t)$
for every $s$,$t\in [0,T]\cap [-h,T-h]$ by uniqueness of the solution.
So we can set
$$G(t-s)=\tilde K(s,t).$$
Replacing in the equation, performing a change of variable
and setting $x=t-s$,
one gets
$$\int_0^x J_0\big(C(x-u)\big)
G(u)~du =-C J_1(Cx) \quad\mbox{ for every }x\in [-T,T].
$$
As $G$ is easily seen to be even this is equivalent to
\begin{equation}\label{equ 28}
\int_0^x J_0\big(C(x-u)\big)
G(u)\,du=-C J_1(Cx)\quad\mbox{for every }x\in [0,T].
\end{equation}
Performing a Laplace transformation in $x$ gives the equality of functions
$${\cal L}J_0(C\cdot)\times {\cal L}G={\cal L}(-CJ_1(C\cdot)),
$$
from which we deduce the Laplace transform ${\cal L}G$ of $G$,
and by the inverse Laplace transform (see for instance the formulas for
Laplace transformations in \cite{zemanian})
we obtain that
$$G(x)=-\frac{CJ_1(Cx)}{x} {\rm\ on\ }[0,T],
$$
hence on $[-T,T]$ by parity,
and eventually one gets
\begin{equation}\label{equ 29}
\tilde K(s,t)=-\frac{CJ_1(C(t-s))}{t-s}
{\rm\ for\ every\ }s,t\in [0,T],
\end{equation}
giving Theorem \ref{Thsuppl2}. \hfill$\Box$ \smallskip

Concerning the growth in terms of power of $k$, we see that
$$|\tilde K(s,t,y,ik)|
\leq |k_\perp|^2 \|E\|_\infty \|J_1'\|_\infty
\leq (1+|k|^2)\|E\|_\infty.$$

Now we shall prove Theorem \ref{Thsuppl1}. Here $E$
can be written in the form
$$E_\perp(t,\tau,x_\perp)={\cal E}(t,x_\perp)g(\tau,x_\perp),$$
with ${\cal E}$ $\mathbb{R}$-valued and $g$ $\mathbb{R}^2$-valued.
For simplicity we fix $x$
and we write
$$E_\perp(t,\tau,x_\perp)=E_\perp(t,\tau)={\cal E} (t)g(\tau).
$$
Then $\tilde E(t)={\cal E} (t)V$, where $V$ is the vector
$V=\frac{1}{2\pi}\int_0^{2\pi} r\big(g(\tau),\tau\big)\,d\tau$,
and
$$ L_\perp(s,t)=V\int_s^t{\cal E} (\sigma)\,d\sigma=V\big(\Phi(t)
-\Phi(s)\big)$$
where we note $\Phi$ a primitive of ${\cal E} $.
The equation (\ref{equ 27}) becomes
\begin{multline*}
 \int_s^t J_0\big(|k_\perp||V||\Phi(\sigma)-\Phi(s)|\big)
\tilde K(\sigma,t,ik)\,d\sigma \\
 =-\frac{\Phi(t)-\Phi(s)}{|\Phi(t)-\Phi(s)|}|k_\perp||V|{\cal E} (t)
J_1\big(|k_\perp||V||\Phi(t)-\Phi(s)|\big),
\end{multline*}
in other words
\begin{multline}\label{equ 30}
 \int_s^t J_0\big(|k_\perp||V|(\Phi(\sigma)-\Phi(s))\big)
\tilde K(\sigma,t,ik)\,d\sigma \\
 =-|k_\perp||V|{\cal E} (t)
J_1\big(|k_\perp||V|(\Phi(t)-\Phi(s))\big)
\end{multline}
for every $s,t\in [0,T]$.
Now if ${\cal E} (t)$ is strictly positive
(respectively strictly negative)
for every $t\in ]a,b[$,
we have $\Phi'(t)={\cal E} (t)>0$ (respectively $<0$)
hence $\Phi$ is bijective from $[a,b]$
onto the interval $I$ of extremities $\Phi(a)$, $\Phi(b)$.
We make the change of variables $u=\Phi(\sigma)$
in the integral
and we obtain
\begin{multline*}
 \int_\alpha^\beta J_0\big(|k_\perp||V|(u-\alpha)\big)
\tilde K(\Phi^{-1}(u),\Phi^{-1}(\beta),ik)~\frac{du}{{\cal E} (\Phi^{-1}(u))}\\
= -|k_\perp||V|{\cal E} (\Phi^{-1}(\beta))
J_1\big(|k_\perp||V|(\beta-\alpha)\big)
{\rm\ for\ every\ }\alpha,\beta\in I,
\end{multline*}
where we have set $\alpha=\Phi(s)$ and $\beta=\Phi(t)$.
We set
$$\Gamma(\alpha,\beta)=\frac{\tilde K(\Phi^{-1}(u),\Phi^{-1}(\beta))}
{{\cal E} (\Phi^{-1}(\alpha)){\cal E} (\Phi^{-1}(\beta))}$$
and then we have
\begin{multline*}
 \int_\alpha^\beta J_0\big(|k_\perp||V|(u-\alpha)\big)
\Gamma(u,\beta)\,du \\
 =-|k_\perp||V|J_1\big(|k_\perp||V|(\beta-\alpha)\big)
{\rm\ for\ every\ }\alpha,\beta\in I.
\end{multline*}
We met the same equation in the preceding case. The solution
is
$$\Gamma(\alpha,\beta)=-\frac{|k_\perp||V|
J_1\big(|k_\perp||V|(\beta-\alpha)\big)}{\beta-\alpha}
,$$
which gives
\begin{equation}\label{equ 31}
\tilde K(s,t)=-{\cal E} (t){\cal E} (s)\frac{|k_\perp||V|
J_1\big(|k_\perp||V|(\Phi(t)-\Phi(s))\big)}{\Phi(t)-\Phi(s)}
{\rm\ for\ every\ }s,t\in [a,b].
\end{equation}
At this point we know $\tilde K(s,t)$ for $s$, $t$ both in an interval
where ${\cal E} $ is strictly positive or strictly negative.
For such $s$, $t$ we have, replacing in (\ref{equ 30}), the equality
\begin{multline*}
 \int_s^t J_0\big(|k_\perp||V|(\Phi(\sigma)-\Phi(s))\big)
{\cal E} (\sigma)\frac{
J_1\big(|k_\perp||V|(\Phi(t)-\Phi(\sigma))\big)}{\Phi(t)-\Phi(\sigma)}
\,d\sigma \\
= J_1\big(|k_\perp||V|(\Phi(t)-\Phi(s))\big).
\end{multline*}
We can write
\begin{multline*}
 J_0\big(|k_\perp||V|(\Phi(\sigma)-\Phi(s))\big)\\
 = \sum_{n=0}^\infty \frac{J_0^{(n)}\big(|k_\perp||V|(\Phi(t)-\Phi(s))\big)}{n!}
(-1)^n |k_\perp|^n|V|^n(\Phi(t)-\Phi(\sigma))^n,
\end{multline*}
$$\frac{J_1\big(|k_\perp||V|(\Phi(t)-\Phi(\sigma))\big)}{\Phi(t)-\Phi(\sigma)}=
\sum_{n=1}^\infty \frac{J_1^{(n)}(0)}{n!}
 |k_\perp|^n|V|^n(\Phi(t)-\Phi(\sigma))^{n-1},$$
and write the product of the two series as
\begin{multline*}
 J_0\big(|k_\perp||V|(\Phi(\sigma)-\Phi(s))\big)
\frac{J_1\big(|k_\perp||V|(\Phi(t)-\Phi(\sigma))\big)}{\Phi(t)-\Phi(\sigma)}\\
 =\sum_{n=0}^\infty P_n(s,t)(\Phi(t)-\Phi(\sigma))^{n}.
\end{multline*}
Integrating term by term we get
$$\int_s^t \sum_{n=0}^\infty P_n(s,t)(\Phi(t)-\Phi(\sigma))^{n}
{\cal E} (\sigma)\,d\sigma =\sum_{n=0}^\infty
P_n(s,t)\frac{(\Phi(t)-\Phi(s))^{n+1}}{n+1},$$ so that the
following equality holds:
$$\sum_{n=0}^\infty P_n(s,t)\frac{(\Phi(t)-\Phi(s))^{n+1}}{n+1}
=\sum_{n=1}^\infty \frac{J_1^{(n)}(0)}{n!}
 |k_\perp|^n|V|^n(\Phi(t)-\Phi(s))^{n}.$$
As this is an equality of series which has nothing to do with the sign of
$\mathcal{E}$ we deduce that the formula (\ref{equ 31}) giving
$\tilde K(s,t)$ is valid for every $s$, $t$ in $[0,T]$.
The inequality
$$|K(s,t,y,ik)|\leq \|{\cal E} \|_\infty |\tilde g(y)|^2(1+|k|^2)$$
is immediate. This completes the proof of Theorem \ref{Thsuppl1}.

\noindent\textsc{Emmanuel Fr\'enod} \\
LMAM, Universit\'e de Bretagne Sud, \\
Campus de Tohannic, F-56000, Vannes, France\\
e-mail:  Emmanuel.Frenod@univ-ubs.fr \smallskip

\noindent\textsc{Fr\'ed\'erique Watbled}\\
LMAM, Universit\'e de Bretagne Sud, \\
Campus de Tohannic, F-56000, Vannes, France \\
and \\
IRMAR, Universit\'e Rennes 1, Campus de Beaulieu,
35042 Rennes cedex\\
e-mail: watbled@maths.univ-rennes1.fr

\end{document}